\documentclass[number,citesort,seceqn,dvips]{arxbj}
\usepackage{graphicx}

\aid{0}
\volume{16}
\issue{4}
\pubyear{2010}
\firstpage{1240}
\lastpage{1261}
\doi{10.3150/09-BEJ242}

\makeatletter
\DeclareMathAlphabet\mathcaligr{OMS}{cmsy}{m}{n}
\renewcommand{\mathcal}{\mathcaligr}
\newcommand{\cal}{\mathcal}

\newproclaim{definition}{Definition}
\newremark{remark}{Remark}
\newtheorem{proposition}{Proposition}
\newtheorem{lemma}{Lemma}
\newremark{APMCM}{Allowance prices via Monte Carlo method}
\newremark{VECMCM}{Valuation of European call via Monte Carlo method}
\newremark{example}{Example}
\makeatother

\begin{document}
\begin{frontmatter}

\title{On fair pricing of emission-related derivatives}
\runtitle{On fair pricing of emission-related derivatives}

\begin{aug}
\author[a]{\fnms{Juri} \snm{Hinz}\corref{}\thanksref{a}\ead[label=e1]{mathj@nus.edu.sg}}
\and
\author[b]{\fnms{Alex} \snm{Novikov}\thanksref{b}\ead[label=e2]{Alex.Novikov@uts.edu.au}}
\runauthor{J. Hinz and A. Novikov}
\address[a]{Department of Mathematics, National University of Singapore,
2 Science Drive, 117543 Singapore. \printead{e1}}
\address[b]{Department of Mathematical Sciences,
University of Technology Sydney, PO Box 123,
Broadway, NSW 2007, Australia. \printead{e2}}
\end{aug}

\received{\smonth{7} \syear{2008}}
\revised{\smonth{9} \syear{2009}}

%
\begin{abstract}
Tackling climate change is at the top of many agendas. In this context,
emission trading schemes
are considered as promising
tools. %
The regulatory framework for an emission trading scheme introduces a
market for emission allowances and creates
a need for risk management by appropriate financial contracts.
In this work, we address logical principles underlying their valuation.
\end{abstract}

%

\end{frontmatter}

\section{Introduction}

The generic principle of an emission trading scheme is based on the
so-called `cap-and-trade' mechanism.
In this framework, an authority allocates
fully tradable
credits among responsible institutions. At pre-settled compliance dates,
each source must have enough allowances
to cover all of its recorded emissions, or be subject to penalties.

A mandatory cap-and-trade system involves its participants
in a risky business with an obvious need for risk management. That is,
certificate trading
is usually accompanied by a secondary market
for emission-related futures, including a rapidly growing variety of their
derivatives.
Their pricing is addressed in this work.

%
Our contribution focuses on a methodology between equilibrium and risk-neutral
approaches. Due to the complexity of emissions markets,
risk-neutral dynamics must be addressed in terms of explanatory variables,
viewed as proxies of fundamental quantities.
Thus, we utilize equilibrium analysis to explain the role of
fundamentals in
risk-neutral allowance price formation. Thereby, the key issue is
a feedback relation between allowance prices and abatement activity.
Namely, we demonstrate that
any increase in allowance
price causes market participants to enforce emission saving in order to
sell their allowances. Hence,
an increasing allowance price encourages a supply of certificates and
lowers the probability of non-compliance, which
tends to bring down their prices.
Apparently, the correct description of this feedback is the key to
derivatives pricing.
The present work focuses on this issue. On this account, our contribution
goes beyond any risk-neutral approach to modeling of emission-related
assets suggested in the existing literature to date.

\section{Emissions markets}

The literature on this subject is enormous: it encompass
hundreds of books and papers.
For this reason, we focus only on those market models which are
relevant in the present approach.

\textit{Economic theory} of allowance trading can be traced back to \cite
{Dales} and \cite{Montgomery},
whose authors proposed a market model for the public  environmental goods described by tradable permits.

\textit{Dynamic allowance trading} is addressed in \cite{Cronshaw,Tietenberg2,Rubin,Leiby,Schennach,StevensRose,Maeda}
and in the literature cited therein.

\textit{Empirical evidence} from existing markets is discussed in
\cite{DaskalakisPsychoyiosMarkellos}. This paper suggests economic
implications and hints at several ways to model spot and futures
allowance prices, whose detailed interrelations are investigated in
\cite{UhrigHomburgWagner1} and
\cite{UhrigHomburgWagner2}. 

\textit{Econometric modeling} is addressed in
\cite{BenzTrueck}, where characteristic properties for financial time
series are observed for prices of emission allowances from
the mandatory European Scheme EU ETS. Furthermore, a Markov switch and AR--GARCH
models are suggested. The work \cite{PaolellaTaschini} also considers
tail behavior and the heteroscedastic dynamics in the returns of
emission allowance prices.

\textit{Dynamic price equilibrium and optimal market design} are
investigated in
\cite{CarmonaFehrHinz}. 
Based on this approach, \cite{CarmonaFehrHinzPorchet} discusses the
price formation
for goods whose production is affected by emission regulations. In this
setting, an equilibrium analysis
confirms the existence of the so-called `windfall profits' (see \cite
{sijm}) and provides
quantitative tools to analyze alternative market designs.

\textit{Pricing of options} was addressed only recently. The paper \cite
{ChesneyTaschini} discusses an endogenous emission
permit price dynamics within equilibrium setting and elaborates on
valuation of European option on emission allowances.
The paper \cite{Seifert} and the dissertation \cite{Wagner} deal with
the the risk-neutral allowance price formation
within EU ETS. Here, utilizing equilibrium properties, the price
evolution is
treated in terms of marginal abatement costs and optimal stochastic control.
Also, the work \cite{CetinVerschuere} is devoted to option pricing
within EU ETS.
The authors suppose that the drift of allowance spot prices is related
to a hidden variable
which describes the overall market position in allowance contracts and
they make use of
filtering techniques to derive option price formulas which reflect
specific allowance
banking regulations, valid in the EU ETS. Finally, the recent work
\cite{CarmonaHinz} presents an approach
where emission certificate futures are modeled in terms of
deterministic time change applied to a certain class of interval-valued
diffusion processes.
The present work brings aspects of risk-aversion into the line of
research followed in \cite{Seifert,CarmonaFehrHinzPorchet}
and \cite{CarmonaFehrHinz}, which we briefly sketch now. %
Within a stochastic model of an emissions market, a so-called
\textit{central planer} problem is introduced and discussed in
\cite{Seifert}. Under additional assumptions, the authors formulate
this problem in terms of continuous-time stochastic
optimization.
Furthermore, they provide economic arguments justifying
why optimal control solutions correspond to an equilibrium of the
emissions market.
Interpreting the allowance certificate price as the marginal abatement
costs, particular explicit solutions
are discussed and yield a dynamic stochastic model for allowance price
evolution.
The work \cite{CarmonaFehrHinz} starts from the opposite direction. In
a discrete-time
framework, the Radner equilibrium of an emissions market is introduced
and constructed
via a solution of the central planer problem.
The work \cite{CarmonaFehrHinzPorchet} yields an extension: in a
slightly different setting,
it is proved that any market equilibrium is reached by this
methodology. 
Thus, results from \cite{Seifert,CarmonaFehrHinz} and \cite
{CarmonaFehrHinzPorchet} show that a
quantitative analysis of emissions markets is tractable in terms of
stochastic control theory. 
However, this connection
is valid only if risk aversion is neglected, in other words, under the
assumption that \textit{each agent} possesses a \textit{linear utility function}.
Losing sight of risk aversion comes at the costs of
unrealistic results. Among other singularities, it turns out that the
equilibrium allowance price follows
a martingale (with respect to objective measure!) with the consequence
that allowance trading can be
arbitrary, only the final position must be adjusted accordingly.

This work resolves all of these problems. Starting from the
no-arbitrage property which is satisfied in an equilibrium of a
market with risk-averse players,
we show that the risk-neutral allowance price dynamics exhibits the
above feedback property,
which we formalize as a fixed point equation, discussing its solution.
We show that for such a risk-averse setting, our fixed point equation
plays the same role as
the central planer optimal control problem for the non-risk-averse
situation. Namely, it provides a
methodology to describe the market equilibrium in terms of \textit{aggregated quantities}.
However, this description is valid only from the viewpoint of the
so-called risk neutral dynamics, not being suitable
for discussing all interesting problems. Still,
derivatives valuation
is naturally addressed in, and can be obtained in, this setting.


\section{Mathematical model}

Let $(\Omega, {\cal F}, {\mathrm P}, ({\cal F}_{t})_{t=0}^{T})$ be a filtered
probability space.
Assume that ${\cal F}_{0}$ is deterministic and agree that all processes
considered in this work are
adapted to $({\cal F}_{t})_{t=0}^{T}$. Write $\mathbb{E}_t(\cdot)$ and $
{\mathrm P}_{t}$
to denote, respectively, conditional expectation and conditional
distribution with respect to ${\cal F}_t$.
Consider a market with a finite number $I$ of the agents confronted
with emission reduction.
%

\textit{Emission dynamics.}
For each agent
$i \in I$, introduce the stochastic process $(E^{i}_{t})_{t=0}^{T-1}$
with the interpretation that
$E^{i}_{t}$ describes the total pollution of the agent $i$ which is
emitted within the time
interval $[t, t+1]$
in the case of the so-called `business-as-usual' scenario (where no
abatement measure is applied).
Although each agent is considered as a potential producer, purely
financial institutions
are also covered with this approach by setting emissions to zero, that
is, $E^{i}_{t}=0$ for $t=0, \ldots, T-1$.

\textit{Abatement.}
Consider the opportunity to reduce emissions. Each agent $i$
can decide at any time $t=0, \ldots, T-1$ to reduce its emissions
within $[t, t+1]$ by $\xi^{i}_{t}$ pollutant units.
We suppose that each abatement level is possible, ranging from no
reduction to full reduction.
Hence, we assume that $0 \le\xi^{i}_{t} \le E^{i}_{t}$ holds for all
$t=0, \ldots, T-1$.

\textit{Abatement costs.} We
assume that the cost of abatement is a random function of the reduced
volume. The randomness
is due to uncertainty in prices (of fuel) and is observable at the
corresponding time. Thus, if
the agent $i$ decides at time $t=0, \ldots, T-1$ on reduction of
their own emissions by $x \in[0, \infty[$ units, then
it causes costs $C^{i}_{t}(x)$, where given
%
\begin{equation}
\cases{
C^{i}_{t} \dvtx [0, \infty[ \times\Omega\mapsto\mathbb R \mbox{ is }{\cal
B}([0, \infty[)\otimes{\cal F}_{t}\mbox{-measurable}\cr
\mbox{and for each }\omega\in\Omega, x \mapsto C^{i}_{t}
(x)(\omega) \mbox{ is strictly}\cr
\mbox{convex and continuous with }C(0)=0.
}
\label{asscosts}
\end{equation}
Since emission savings cannot exceed the business-as-usual emission,
the abatement activity $(\xi^{i}_{t})_{t=0}^{T-1}$ is feasible if
%
\begin{equation} \label{mustsatisfy}
0 \le\xi^{i}_{t} \le E^{i}_{t}, \qquad t =0, \ldots, T-1.
\end{equation}
Following abatement policy $(\xi^{i}_{t})_{t=0}^{T-1}$, the agent $i$
accumulates
at the compliance date $T$ the total terminal costs
%
\begin{equation}
\sum_{t=0}^{T-1}C^{i}_{t}(\xi^{i}_{t}). \label{abatement}
\end{equation}

\textit{Abatement volume.}
For later use, let us introduce, for each $\omega\in\Omega$, $t=0,
\ldots, T-1$ and $a \in[0, \infty[$,
the abatement volume $c_{t}^{i}(a)(\omega)$ as
%
\begin{equation} \label{volumefunction}
c_{t}^{i}(a)(\omega):=\arg\min\{ C_{t}^{i}(x)(\omega)-ax  \dvtx
x \in[0, E^{i}_{t}(\omega)]\},
\end{equation}
which is well defined since, under the assumptions (\ref{asscosts}),
the minimum of the function $x \mapsto C_{t}^{i}(x)(\omega)-ax$ on
$[0, E^{i}_{t}(\omega)]$ is attained at the unique point.
The reader may imagine $c_{t}^{i}(a)(\omega)$ as 
the total reduction volume which is available within $[t, t+1]$ in the
situation $\omega$
at a price which is less than or equal to $a$ (measured in currency
unit per pollutant unit).
A straightforward proof shows that (\ref{asscosts}) ensures that
%
\begin{eqnarray}
&&[0, \infty[ \mapsto\mathbb R,  a \mapsto c_{t}^{i}(a)(\omega)
\mbox{ is non-decreasing and}\nonumber\\[-11pt]\\[-11pt]
&&\mbox{continuous for almost every }\omega\in\Omega \mbox{ and }t=0, \ldots
, T-1.\nonumber
\label{nondecrcont}
\end{eqnarray}
%
For later use, we
introduce the cumulative abatement volume function
%
\begin{equation} \label{cfunt}
c_{t}(a):=\sum_{i \in I}c^{i}_{t}(a),\qquad  a \in[0, \infty[.
\end{equation}
Obviously, $c_{t}(a)(\omega)$ stands for the total abatement in the
market, which is available from
all measures in the situation $\omega$ whose price is less than or
equal to $a \in[0, \infty[$.

\textit{Allowance trading.}
Suppose that, at any time $t=0, \ldots, T$, credits can be
exchanged between agents by trading at the spot price $A_{t}$.
Denote by $\vartheta^{i}_{t}$ the change at time $t$ in allowance
number held
by agent $i$. That is, given the allowance prices $(A_{t})_{t=0}^{T}$,
the position changes $(\vartheta^{i}_{t})_{t=0}^{T}$ yield costs
%
\begin{equation}
\sum_{t=0}^{T}\vartheta^{i}_{t}A_{t}. \label{trading}
\end{equation}

\textit{Penalty payment.}
The total pollution of the agent $i$ can be expressed as a difference
\[
\sum_{t=0}^{T-1}E^{i}_{t} - \sum_{t=0}^{T-1}\xi^{i}_{t}
\]
of the cumulative business-as-usual emission less the entire reduction.
As mentioned above,
a penalty $\pi\in]0, \infty[$ is being paid at maturity $T$ for each
unit of pollutant, which is
not covered by allowances.
Considering the total change in the allowance position $\sum
_{t=0}^{T}\vartheta^{i}_{t}$ effected by trading, the
loss of the agent $i$ resulting from potential penalty payment is
%
\begin{eqnarray}
\pi\Biggl(\sum_{t=0}^{T-1}(E^{i}_{t}-\xi^{i}_{t}-\vartheta
^{i}_{t})-\gamma^{i}- \vartheta^{i}_{T}\Biggr)^{+}, \label{penaltypayment}
\end{eqnarray}
where
%
\begin{equation}
\gamma^{i},\ i \in I\mbox{ are agents' initial allowance
allocations}. \label{initall}
\end{equation}
\begin{remark}
Our stylized scheme deals with stand-alone emission trading mechanisms.
In the real world,
cap-and-trade systems operate on multi-period scales, where unused allowances
can be carried out (banked) into next period. Further period
interconnections may include
a transfer of future allocation from the next into the present period
(borrowing) and,
in the case of non-compliance, a withdrawal of an appropriate number of
credits from the
next period allocation in addition to penalty payment. To complete the complexity,
let us mention that different emissions markets could be interconnected
by acceptance of foreign certificates in the national scheme.
Emission trading in multi-period settings is addressed in, among
others, \cite{CarmonaHinz} and
\cite{CetinVerschuere}. Mathematically, it reduces to the
specification of a more complex
penalty mechanism than that presented above. For this reason, we have
decided to focus on the stand-alone
allowance market to analyze quantitative methods in the simplest situation
before tackling multi-scale systems (such as the second period of EU ETS).
\end{remark}

\textit{Recording uncertainty.}
In what follows, we also need to take into account uncertainty in the
emission recording.
It is convenient to subtract these recording errors from the initial
allocation. Hence, we
interpret $\gamma^{i}$ as the credits allocated to the agent $i$ less
emissions which
become known with certainty only at time $T$. With this interpretation,
$\gamma_{i}$ stands
for allowances effectively available for compliance and is modeled by
an ${\cal F}_{T}$-measurable random variable. For later use, let us
agree that
the distribution of
$\sum_{i\in I}\gamma^{i}$, conditioned on
${\cal F}_{T-1}$, possesses almost surely no point masses,
which implies that
%
\begin{equation}
P\biggl(\sum_{i \in I}\gamma^{i}=X\biggr)=0\qquad \mbox{for each ${\cal F}
_{T-1}$-measurable $X$.} \label{nopointmassnew}
\end{equation}

\textit{Admissible policies.}
Since maximally possible reduction cannot exceed emission, we\break have~(\ref{mustsatisfy}).
Let us define the space of feasible trading $\vartheta^{i}=(\vartheta
^{i}_{t})_{t=0}^{T}$
and abatement strategies $\xi^{i}=(\xi^{i}_{t})_{t=0}^{T-1}$ of the
agent $i \in I$ by
%
\begin{eqnarray}
{\cal U}^{i}:=\{ (\vartheta^{i}, \xi^{i})   \dvtx
0 \le\xi^{i}_{t} \le E^{i}_{t},  t =0, \ldots, T-1\}. \label{admiss}
\end{eqnarray}

\textit{Individual wealth.}
In view of
(\ref{abatement}), (\ref{trading}) and (\ref{penaltypayment}), the
revenue of the agent $i$ following admissible
policy $(\vartheta^{i}, \xi^{i}) \in{\cal U}^{i}$ equals
%
\begin{eqnarray} \label{totrev}
L^{A,i}(\vartheta^{i}, \xi^{i})&=&-\sum_{t=0}^{T-1}\bigl(\vartheta
^{i}_{t}A_{t}+ C^{i}(\xi^{i}_{t})\bigr)\nonumber\\[-8pt]\\[-8pt]
&&{} -\vartheta^{i}_{T}A_{T}
- \pi\Biggl(\sum_{t=0}^{T-1}(E^{i}_{t}-\xi^{i}_{t}-\vartheta^{i}_{t}) -
\gamma^{i}-\vartheta^{i}_{T}\Biggr)^{+}.\nonumber
\end{eqnarray}

\textit{Risk aversion.} To face risk preferences, suppose that attitudes
of the agents $i \in I$ are
described by utility functions
$
U^{i} \dvtx \mathbb R\mapsto\mathbb R,
$ which are continuous, strictly increasing and concave.
Consider the utility functional
$
u^{i}(X)=E(U^{i}(X)),
$
which is assumed to be defined for each random variable $X$ where the
expectation is finite or $+ \infty$.
Given allowance price process $A=(A_{t})_{t=0}^{T}$, the agent $i$
behaves rationally,
maximizing
$
(\vartheta^{i}, \xi^{i}) \mapsto u^{i}(L^{A,i}(\vartheta^{i}, \xi
^{i}))$ by an appropriate choice
of their own policy $(\vartheta^{i*}, \xi^{i*})$.

\textit{Market equilibrium.}
Following standard theory, a realistic market state is described by the
so-called equilibrium --
a situation where the allowance price, positions and abatement measures
are such that each agent is satisfied by their own policy and,
at the same time, natural restrictions are fulfilled. In our framework,
an appropriate notion of equilibrium
is given as follows.

\begin{definition}
\label{de:equilibrium}
The process $A^{*}=(A^{*}_{t})_{t=0}^{T}$
is called an \textit{equilibrium allowance price process} if, for each
$i \in I$,
there exists $(\vartheta^{*i}, \xi^{*i}) \in{\cal U}^{i} $
such that $u^{i}(L^{A^{*}, i} (\vartheta^{*i}, \xi^{*i}))$ is finite
and
\begin{enumerate}[(ii)]
\item[(i)] the cumulative changes in positions are in zero net supply,
that is,
%
\begin{equation} \label{eqold}
\sum_{i \in I}\vartheta_{t}^{*i}=0 \qquad \mbox{for all $t=0, \ldots, T$};
\end{equation}
\item[(ii)] each agent $i \in I $ is satisfied by their own policy, in
the sense that
%
\begin{eqnarray}\label{IndividualOptimum}
&&u^{i}(L^{A^{*}, i} (\vartheta^{*i}, \xi^{*i})  )\ge
u^{i}(L^{A^{*}, i}
(\vartheta^{i}, \xi^{i})
 ) \nonumber\\[-10pt]\\[-10pt]
&&\quad\mbox{for each }(\vartheta^{i}, \xi^{i}) \in{\cal U}^{i}\mbox{ where }
u^{i}(L^{A^{*}, i}
(\vartheta^{i}, \xi^{i})
 )\mbox{ exists.}\nonumber
\end{eqnarray}
\end{enumerate}
\end{definition}

The existence of emissions market equilibrium is addressed in \cite
{CarmonaFehrHinz} and \cite{CarmonaFehrHinzPorchet},
under the assumption of a linear utility function and in a slightly
different setting.
However, although equilibrium modeling in the spirit of these
contributions is appropriate
to investigate important questions of
optimal market design, it
has little to offer to the problem of derivatives valuation.
With the present approach, we intend to establish a reduced-form model
which describes the evolution of emission-related assets
from a risk-neutral perspective.
We obtain a realistic picture
by incorporating three essential assumptions into a risk-neutral model.
These assumptions are shown to be direct consequences of an equilibrium
situation:
\begin{enumerate}[(b)]
\item[(a)]
There is no arbitrage since, in equilibrium, any profitable strategy
would immediately be followed by all agents.
This would instantaneously
change prices and exhaust any arbitrage opportunity.
\item[(b)]
The allowance trading instantaneously triggers all abatement measures
whose costs are below allowance price. The explanation here is that if
an agent possess a technology
with lower reduction costs than the present allowance price, then it is
optimal for that agent to immediately reduce
pollution and take profit from selling allowances.
\item[(c)]
There are only two final outcomes for allowance price. Either the
terminal allowance price drops to zero
or it approaches the penalty level.
The reason is that at maturity, the price must vanish if there is an
excess in allowances, whereas in the
case of their shortage, the price will rise, reaching penalty.
We believe that in reality, an exact coincidence of allowance demand
and supply
occurs with zero probability and can be neglected.
\end{enumerate}
Let us formalize the above assertions (a), (b) and (c).
\begin{proposition} \label{mainprop}
Suppose that $(A^{*}_{t})_{t=0}^{T}$ is an equilibrium allowance price
and $(\xi^{i*}_{t})_{t=0}^{T-1}$ for $i \in I$ are the corresponding
equilibrium abatement policies.
\begin{enumerate}[(b)]
\item[(a)] There exists
a measure $Q$ which is equivalent to $P$ such that
$(A^{*}_{t})_{t=0}^{T}$ follows a $Q$-martingale.
\item[(b)] For each $i \in I$, we have
%
\begin{equation}
\xi^{i*}_{t}= c_{t}^{i}(A^{*}_{t}), \qquad t=0, \ldots, T-1,
\label{marginalcosts}
\end{equation}
with abatement volume functions $c_{t}^{i}$, $t=0, \ldots, T-1$, from
(\ref{volumefunction}).
\item[(c)] The terminal value of the allowance price is given by
%
\begin{equation} \label{termialprice}
A^{*}_{T}=\pi1_{\{ \sum_{i \in I}(\sum_{t=0}^{T-1}(E^{i}_{t}-\xi
^{i*}_{t})- \gamma^{i}) \ge0 \}}.
\end{equation}
\end{enumerate}
\end{proposition}

Before we proceed with the proof, let us emphasize that this result can
serve as a starting point for risk-neutral modeling. The above
proposition states that at equilibrium, the allowance price process
$(A^{*}_{t})_{t=0}^{T}$
follows a martingale with respect to an equivalent measure $ {\mathrm
Q}\sim {\mathrm P}
$ whose terminal value is
\[
A^{*}_{T}=\pi1_{\{ \sum_{i \in I}\sum_{t=0}^{T}(E^{i}_{t}-\xi
^{*}_{t}) -\gamma^{i}) \ge0\}},
\]
obviously depending on intermediate values $(A^{*}_{t})_{t=0}^{T-1}$ through
abatement volume function $\xi^{i*}_{t}=c^{i}_{t}(A^{*}_{t})$ for
$t=0, \ldots, T-1,
i \in I$.
The surprising and far-reaching consequence is that, from a
risk-neutral perspective, only cumulative
market quantities are relevant. To see this, define the overall
allowance shortage
%
\begin{eqnarray} \label{capadjusted}
{\cal E}_{T}=\sum_{i \in I}\Biggl(\sum_{t=0}^{T-1}E^{i}_{t} -\gamma^{i}\Biggr)
\end{eqnarray}
which would appear in the market without any emission penalty.
Further, recall from (\ref{volumefunction})
and~(\ref{cfunt}) the cumulative abatement functions
to express the risk-neutral certificate price dynamics
in terms of the following feedback equation:
\[
A_{t}=\mathbb{E}^{Q}_{t}\bigl(\pi1_{\{{\cal E}_{T}-\sum_{t=0}^{T-1}c_{t}(A^{*}_{t})
\ge
0\}} \bigr), \qquad t=0, \ldots, T-1.
\]
Although individual market attributes and actions of the different
agents seem to be irrelevant in this picture,
the reader should notice that this picture appears only from the
risk-neutral viewpoint. In line with standard aggregation theorems, the
equilibrium market state heavily depends on, and is determined by,
market architecture, rules, risk attitudes
and uncertainty. However, once equilibrium is reached and all arbitrage
opportunities are
exhausted, asset dynamics can be considered under risk-neutral measure.
With respect to this measure, market evolution appears as if it were driven
by cumulative quantities only.

With this in mind, let us formulate the problem of the
reduced-form modeling as
follows:
%
\begin{equation}
\cases{
\mbox{Given measure } {\mathrm Q}\sim {\mathrm P},\mbox{ random variable }
{\cal E}_{T} \cr
\mbox{and abatement volume functions }(c_{t})_{t=0}^{T-1},\cr
\mbox{determine a } {\mathrm Q}\mbox{-martingale } (A^{*}_{t})_{t=0}^{T}
\mbox{ with} \cr
A^{*}_{T}=\pi1_{\{{\cal E}_{T}-\sum_{t=0}^{T-1}c_{t}(A^{*}_{t})
\ge0\}}.
}\label{ReducedProblem}
\end{equation}
Note that this formulation serves as a guideline for martingale modeling
since price-dependent abatement volume $c_{t}(a)$ can be estimated from
market data,
whereas potential allowance shortage ${\cal E}_{T}$
can be modeled in terms of total allowance allocation and demand
fluctuations on goods whose
production causes the pollution. Finally, we shall emphasize a natural
passage to continuous time.
%
\begin{equation}
\cases{
\mbox{Given, on a probability space }(\Omega, {\cal F}, P, ({\cal F}_{t})_{t
\in[0, T]}), \cr
\mbox{an equivalent measure } {\mathrm Q}\sim {\mathrm P},\mbox{ random
variable }{\cal E}_{T}\cr
\mbox{and a family of abatement functions }(c_{t})_{t \in[0, T]}, \cr
\mbox{determine a } {\mathrm Q}\mbox{-martingale }(A^{*}_{t})_{t \in[0,
T]}\mbox{ with}\cr
A^{*}_{T}=\pi1_{\{{\cal E}_{T}-\int_{0}^{T}c_{t}(A^{*}_{t})\,\mathrm{d}t
\ge 0\}}.
}\label{ReducedProblemCont}
\end{equation}
%
%
\begin{pf*}{Proof of Proposition~\ref{mainprop}}
(a) According to the first fundamental theorem of asset pricing (see
\cite{Kabanov}), it suffices to verify that if $(A^{*}_{t})_{t=0}^{T}$
is an equilibrium allowance price process, then there is no arbitrage
for allowance trading. Let us follow an
indirect proof, supposing that $(\nu_{t})_{t=0}^{T-1}$ is an allowance
trading arbitrage, meaning that
%
\begin{equation}
{\mathrm P}\Biggl(\sum_{t=0}^{T-1}\nu_{t}(A_{t+1}-A_{t}) \ge0\Biggr)=1,\qquad
 {\mathrm P}\Biggl(\sum
_{t=0}^{T-1}\nu_{t}(A_{t+1}-A_{t})>0\Biggr)>0. \label{arbitrage}
\end{equation}
Now, we verify that in the presence of arbitrage, no equilibrium can
exist since each agent $i$ can change their own policy $(\vartheta
^{*i}, \xi^{*i})$
to an improved strategy $(\tilde\vartheta^{i}, \xi^{*i})$ satisfying
%
\begin{equation} \label{utilityincrease}
u^{i}(L^{A^{*}, i} (\vartheta^{*i}, \xi^{*i})  )<
u^{i}(L^{A^{*}, i}
(\tilde\vartheta^{i}, \xi^{*i})
 ).
\end{equation}
The improvement is achieved by incorporating arbitrage $(\nu
_{t})_{t=0}^{T-1}$ into their own allowance trading as follows:
\[
\tilde\vartheta^{i}_{t}:=\vartheta^{*i}_{t}+ (\nu^{i}_{t}-\nu
^{i}_{t-1})\qquad \mbox{for all $t=0,
\ldots, T$,}
\]
with appropriate definitions $\nu_{-1}=\nu_{T}:=0$. Indeed, the
revenue improvement from allowance trading is
\[
-\sum_{t=0}^{T}\tilde\vartheta^{i}_{t}A_{t}=-\sum_{t=0}^{T}
\vartheta^{*i}_{t} A_{t} + \sum_{t=0}^{T-1}\nu
^{i}_{t}(A_{t+1}-A_{t}),
\]
which we combine with (\ref{arbitrage}) to see that there is no
optimality since
\[
{\mathrm P} \bigl(L^{A,i}(\vartheta^{*i}, \xi^{i}) \le L^{A,i}(\tilde
\vartheta^{i}, \xi^{i})  \bigr)=1,\qquad
{\mathrm P} \bigl(L^{A,i}(\vartheta^{*i}, \xi^{i})<L^{A,i}(\tilde
\vartheta
^{i}, \xi^{i})  \bigr)>0
\]
together imply that
\[
u^{i}(L^{A,i}(\vartheta^{*i}, \xi^{i}))<u^{i}(L^{A,i}(\tilde
\vartheta^{i}, \xi^{i})).
\]

(b) To prove (\ref{marginalcosts}), consider the bijection
%
\begin{equation}
{\cal U}^{i} \to{\cal U}^{i},\qquad (\theta^{i}, \xi^{i}) \mapsto
(\phi(\theta^{i},\xi^{i}), \xi^{i}),
\end{equation}
where the transformed trading strategy $\vartheta^{i}=\phi(\theta
^{i},\xi^{i})$ is given by
\[
\vartheta^{i}_{t}=\theta^{i}_{t}-\xi^{i}_{t},\qquad  t=1, \ldots,
T-1,\qquad \vartheta^{i}_{T}=\theta^{i}_{T}.
\]
Obviously, $(\vartheta^{i*}, \xi^{i*})$ is a maximizer to the
original problem
\[
{\cal U}^{i} \to\mathbb R,\qquad (\vartheta^{i}, \xi^{i}) \mapsto
u^{i}(L^{A^{*}, i}(\vartheta^{i}, \xi^{i}) )
\]
if and only if
$(\vartheta^{i*}, \xi^{i*})= (\phi(\theta^{i*},\xi^{i*}), \xi
^{i*})$, where $(\theta^{i*}, \xi^{i*})$ is a
maximizer to the transformed problem
%
\begin{equation} \label{reparam} {\cal U}^{i} \to\mathbb R,\qquad
(\theta^{i}, \xi^{i}) \mapsto u^{i}(L^{A^{*}, i}(\phi(\theta
^{i},\xi^{i}), \xi^{i}) ).
\end{equation}
The last line in the calculation
%
\begin{eqnarray}
&&\hspace*{-25pt}L^{A^{*}, i}(\phi(\theta^{i},\xi^{i}), \xi^{i})
\nonumber\\
&&\hspace*{-25pt}\quad=
-\sum_{t=0}^{T-1}(\theta^{i}_{t}-\xi^{i}_{t})A^{*}_{t}-\sum
_{t=0}^{T-1}C^{i}_{t}(\xi^{i}_{t})
- \pi\Biggl(\sum_{t=0}^{T-1}\bigl(E^{i}_{t} -\xi^{i}_{t}-(\theta^{i}_{t}-\xi
^{i}_{t})\bigr) - \gamma^{i}-\theta^{i}_{T}\Biggr)^{+} \nonumber\\
&&\hspace*{-25pt}\quad=-\sum_{t=0}^{T-1}\theta^{i}_{t}A^{*}_{t}
- \pi\Biggl(\sum_{t=0}^{T-1}(E^{i}_{t} -\theta^{i}_{t}) - \gamma
^{i}-\theta^{i}_{T}\Biggr)^{+}
- \sum_{t=0}^{T-1}\bigl(C^{i}(\xi^{i}_{t})-A^{*}_{t}\xi^{i}_{t}\bigr) \label{lastline}
\end{eqnarray}
shows that if $(\theta^{i*}, \xi^{i*})$ is a maximizer to (\ref
{reparam}), then $\xi^{*}$ must satisfy
$\xi^{i*}_{t}:=c^{i}_{t}(A^{*}_{t})$ for $t=0, \ldots, T-1$,
which proves (\ref{marginalcosts}).

(c) This assertion is proved by an argument identical to that given in
\cite{CarmonaFehrHinz}.
\end{pf*}

\section{Reduced-form modeling}
%
In what follows, we propose a solution to the problem of risk-neutral
allowance price modeling~(\ref{ReducedProblem}).
Below, we prove that under the assumptions given above ((\ref
{nopointmassnew}), in particular, is essential),
the problem (\ref{ReducedProblem}) possess a solution. Moreover, we
show how to obtain the
martingale $(A^{*}_{t})_{t=0}^{T}$.

It turns out that
the martingale closed by ${\cal E}_{T}$ plays a crucial role, so
we introduce
\[
{\cal E}_{t}=\mathbb{E}^{Q}({\cal E}_{T}  |  {\cal F}_{t}),\qquad  t=0,
\ldots, T.
\]
For later use, let us also define its increments as
\[
\varepsilon_{t}= {\cal E}_{t}-{\cal E}_{t-1},\qquad t=1, \ldots, T.
\]
Following the intuition that the equilibrium allowance price should be
uniquely determined by the
present time and the general market situation, we express a candidate
for allowance price as
%
\begin{equation} \label{AfuncForm}
A^{*}_{t}(\omega)=\alpha_{t}(G_{t}(\omega))(\omega), \qquad\omega
\in\Omega,\  t=0, \ldots, T,
\end{equation}
with hypothetic functionals
%
\begin{equation} \label{functionals}
\alpha_{t} \dvtx \mathbb R\times\Omega\to[0, \pi], \qquad{\cal
B}(\mathbb R
) \otimes{\cal F}_{t}\mbox{-measurable for }t=0, \ldots T,
\end{equation}
applied to
%
\begin{equation} \label{Demand}
G_{t}={\cal E}_{t}-\sum_{s=0}^{t-1}c_{s}(A^{*}_{s}), \qquad t=0, \ldots
, T.
\end{equation}
According to (\ref{ReducedProblem}),
this approach yields an obvious definition for $\alpha_{T}$:
%
\begin{equation} \label{finalfunc}
\alpha_{T}(g)(\omega)=\pi1_{[0, \infty[}(g),\qquad \omega\in
\Omega,\  g \in\mathbb R.
\end{equation}
Note that, given functionals (\ref{functionals}),
the price process $(A^{*}_{t})_{t=0}^{T}$ is indeed well defined by recursive
application of (\ref{Demand}) and (\ref{AfuncForm}):
%
\begin{eqnarray}
A^{*}_{t}(\omega)&:=&\alpha_{t}(G_{t}(\omega))(\omega), \label
{r1}\\
G_{t+1}(\omega)&:=&G_{t}(\omega)-c_{t}(A^{*}_{t}(\omega))(\omega
)+\varepsilon_{t+1}(\omega), \label{r2} \\
&&\mbox{started at }G_{0}:={\cal E}_{0}. \label{r3}
\end{eqnarray}
Generated by this recursion, the process $(A^{*}_{t})_{t=0}^{T}$
follows a martingale if, for all $t=0, \ldots T-1$ and almost
all $\omega\in\Omega$, the following holds:
\[
\alpha_{t}(g)(\omega)= \mathbb E^{ {\mathrm Q}}_{t}\bigl(\alpha_{t+1}\bigl(g -
c_{t}(\alpha
_{t}(g)(\omega))(\omega)+\varepsilon_{t+1}\bigr)\bigr)(\omega) \qquad\mbox
{for all }g \in\mathbb R,\
\omega\in\Omega.
\]
Indeed, we have
\begin{eqnarray*}
\mathbb E^{ {\mathrm Q}}_{t}(A^{*}_{t+1})(\omega)
&=&
\mathbb E^{ {\mathrm Q}}_{t}\bigl(\alpha_{t+1}(G_{t+1})\bigr)(\omega) \\
&=& \int_{\Omega} \alpha_{t+1}\bigl(G_{t}(\omega
')-c_{t}(A^{*}_{t}(\omega'))(\omega')
+\varepsilon_{t+1}(\omega') \bigr)(\omega') {\mathrm Q}_{t}(\mathrm{d}\omega
')(\omega) \\
&=& \int_{\Omega} \alpha_{t+1}\bigl(G_{t}(\omega)-c_{t}(A^{*}_{t}(\omega
))(\omega)+\varepsilon_{t+1}(\omega')\bigr)(\omega') {\mathrm
Q}_{t}(\mathrm{d}\omega
')(\omega) \\
&=& \mathbb E^{ {\mathrm Q}}_{t}\bigl(\alpha
_{t+1}\bigl(G_{t}(\omega
)-c_{t}(\alpha_{t}(G_{t}(\omega))(\omega))(\omega)+\varepsilon
_{t+1} \bigr)\bigr)(\omega) \\
&=&
\alpha_{t}(G_{t}(\omega))(\omega)=A^{*}_{t}(\omega).
\end{eqnarray*}
In other words, it is sufficient to ensure that
%
\begin{eqnarray}\label{recsolution}
&&\mbox{for each }g \in\mathbb R,\ \alpha_{t}(g)(\omega)\mbox{ solves} \nonumber\\
&& a= \mathbb E^{ {\mathrm Q}}_{t}\bigl(\alpha_{t+1}\bigl(g -
c(a)+\varepsilon_{t+1}\bigr)\bigr)(\omega) \\
&& \mbox{for almost all }\omega\in\Omega.\nonumber
\end{eqnarray}
In the remainder of this section, we will show
that the functionals (\ref{functionals}) are recursively obtained as
the unique solution to (\ref{recsolution}), starting
with $\alpha_{T}$ from (\ref{finalfunc}). First, let us prepare an
auxiliary result dealing with the solution to (\ref{recsolution})
where no conditional information needs to be considered.
\begin{lemma}\label{lema1}
Given
%
\begin{eqnarray}
 c\dvtx \mathbb R&\to&\mathbb R,\qquad \mbox{non-decreasing,
continuous,}
\label{Cfunct} \\
 \alpha_{1} \dvtx \mathbb R\times\Omega&\to&[0, \pi],\qquad
{\cal
B}(\mathbb R) \otimes{\cal F}\mbox{-measurable,} \\
 g &\mapsto&\alpha_{1}(g)(\omega),\qquad
\mbox{non-decreasing for almost all }\omega\in\Omega,
\label{alphafunct1}
\end{eqnarray}
suppose that the random variable $\varepsilon$ satisfies
%
\begin{eqnarray}\label{measuremu}
&&\mathbb R\to[0,\pi], \qquad x \mapsto \mathbb E^{Q}\bigl(\alpha
_{1}(x+\varepsilon
)\bigr)=\int_{\Omega}\alpha_{1}\bigl(x+\varepsilon(\omega')\bigr)(\omega')\mathrm{Q}
(\mathrm{d}\omega') \nonumber\\[-9pt]
\\[-9pt]
&&\phantom{\mathbb R\to[0,\pi], \qquad\,}\quad\mbox{is continuous.}\nonumber
\end{eqnarray}
For each $g \in\mathbb R$, introduce the function $f^{g} \dvtx [0, \pi]
\to\mathbb R
$ given by
%
\begin{eqnarray}\label{ffunction}
f^{g}(a)
&:=&
a- \mathbb E^{Q}\bigl(\alpha_{1}\bigl(g-c(a) +\varepsilon\bigr)\bigr) \nonumber\\[-8pt]\\[-8pt]
&&\hspace*{-10.8pt}=
a-\int_{\Omega}\alpha_{1}\bigl(g-c(a)+\varepsilon(\omega')\bigr)(\omega
'){\mathrm Q}(\mathrm{d}\omega'),\qquad  a \in[0, \pi].\nonumber
\end{eqnarray}
The following assertions then hold:
%
%
\begin{enumerate}[(iii)]
\item[(i)] for each $g \in\mathbb R$, there exists a unique $\alpha
_{0}(g) \in[0, \pi]$ with $f^{g}(\alpha_{0}(g))=0$;
\item[(ii)] the root $\alpha_{0}(g)$ of $f^{g}$ is obtained as a
limit $\alpha_{0}(g)=\lim_{n \to\infty}a^{g}_{n}$ in the standard
bisection method
%
\begin{equation}
a^{g}_{n}=\tfrac{1}{2}(\overline{a}^{g}_{n} + \underline{a}^{g}_{n}),\qquad
\begin{array}{lll} \overline{a}^{g}_{n+1}=a^{g}_{n}, &\qquad\underline
{a}^{g}_{n+1}:=\underline{a}^{g}_{n}, &\qquad
\mbox{if }f^{g}(a^{g}_{n}) \ge0,\vspace*{2pt}\\
\overline{a}^{g}_{n+1}=\overline{a}^{g}_{n}, &\qquad\underline
{a}^{g}_{n+1}:=a^{g}_{n}, &\qquad \mbox{if }f^{g}(a^{g}_{n}) < 0,
\end{array}
\end{equation}
started at $\underline{a}^{g}_{0}:=0$, $\overline{a}^{g}_{0}:=\pi$;
\item[(iii)] the mapping $\mathbb R\to[0, \pi]$, $g \mapsto\alpha
_{0}(g)$ is non-decreasing and continuous.
\end{enumerate}
\label{firstlemma}
\end{lemma}
\begin{pf}
(i) For each $g \in\mathbb R_{+}$, the function $f^{g}$
is continuous due to (\ref{measuremu}) and the continuity (\ref
{Cfunct}) of $c$.
Thus, the existence of a root follows from
the intermediate value theorem because of
%
\begin{equation}
f^{g}(0) \le0, \qquad f^{g}(\pi) \ge0. \label{initialbisection}
\end{equation}
The uniqueness of the root is ensured by the strict monotonic increase
of $f^{g}$. To verify this, observe that
(\ref{alphafunct1}) and (\ref{Cfunct}) imply that
the subtrahend
\[
a \mapsto\int_{\Omega}\alpha_{1}\bigl(g-c(a)+\varepsilon(\omega
')\bigr)(\omega'){\mathrm Q}(\mathrm{d}\omega')
\]
in (\ref{ffunction})
is non-increasing, whereas the minuend $a \mapsto a $ is strictly
increasing in $a$.
{\smallskipamount=0pt
\begin{longlist}
\item[(ii)] The bisection algorithm is properly initialized because of (\ref
{initialbisection}).
Standard arguments ensure its convergence to the root.
\item[(iii)] To show the monotonic increase of $g \mapsto\alpha_{0}(g)$,
suppose that $g'<g$. Then (\ref{alphafunct1}) ensures that for each $a
\in[0, \pi]$,
\[
\int_{\Omega}\alpha_{1}\bigl(g'-c(a)+\varepsilon(\omega')\bigr)(\omega
'){\mathrm Q}(\mathrm{d}\omega') \le\int_{\Omega}\alpha_{1}\bigl(g-c(a)+
\varepsilon(\omega')\bigr)(\omega'){\mathrm Q}(\mathrm{d}\omega'),
\]
giving $f^{g'}(a) \ge f^{g}(a)$ for all $a \in[0, \pi]$, which
implies that $\alpha_{0}(g') \le\alpha_{0}(g)$.\\
\end{longlist}}
\noindent
Now, let us turn to the continuity.
If $\alpha_{0}(g) \in[0, \pi[$, then there exists $\delta>0$ with
$\alpha_{0}(g)+\delta\le\pi$. Due to the strict monotonic increase of
$f^{g}$, we obtain $0<f^{g}(\alpha_{0}(g)+\delta)$.
If $(g_{n})_{n \in\mathbb N} \subset\mathbb R_{+}$ is a sequence with
$\lim_{n
\to\infty}g_{n}=g$, then according to
(\ref{measuremu}),
%
\begin{equation}
\lim_{n \to\infty}f^{g_{n}}\bigl(\alpha_{0}(g)+\delta\bigr) =f^{g}\bigl(\alpha
_{0}(g)+\delta\bigr)>0.
\end{equation}
Hence, there exists $N \in\mathbb N$ such that $f^{g_{n}}(\alpha
_{0}(g)+\delta)>0$ holds for all $n \ge N$.
Thus, we obtain
%
\begin{equation} \label{supconclusion}
\alpha_{0}(g) \in[0, \pi[ \quad\Longrightarrow\quad\limsup_{n \to
\infty} \alpha_{0}(g_{n}) \le\alpha_{0}(g)+\delta,\qquad
\mbox{if }\alpha_{0}(g)+\delta\le\pi.
\end{equation}
Since $\delta>0$ is arbitrarily small and $0 \le\alpha_{0}(g) \le
\pi$, due to (i),
this implication shows that $\alpha_{0}( \cdot)$ is continuous on
each point $g$ with
$\alpha_{0}(g)=0$.
A similar argument yields
%
\begin{equation} \label{infconclusion}
\alpha_{0}(g) \in]0, \pi]\quad \Longrightarrow\quad\liminf_{n \to
\infty} \alpha_{0}(g_{n}) \ge\alpha_{0}(g)-\delta,\qquad
\mbox{if }\alpha_{0}(g)-\delta\ge0.
\end{equation}
Again, since $\delta>0$ is arbitrary, we obtain the continuity of
$\alpha_{0}( \cdot)$ on each point $g$ with
$\alpha_{0}(g)=\pi$.
If $\alpha_{0}(g) \in]0, \pi[$, then the continuity of $\alpha
_{0}(\cdot)$ on $g$ follows by the combination of (\ref
{supconclusion}) and~(\ref{infconclusion}).
\end{pf}

Let us now turn to the conditioned version of Lemma~\ref{lema1}. Supposing the
existence of the regular ${\cal F}_{t}$-conditioned distribution
$Q_{t}$,
the proof reproduces the arguments of the previous lemma with
appropriate notational changes due to
conditioning on
the event $\omega\in\Omega$. However, a useful insight is that the
approximating points $a^{g}_{n}$, $n=0, 1, 2, \ldots,$
of the bisection algorithm turn out to be dependent on $g \in\mathbb
R$ and
$\omega\in\Omega$ in a ${\cal B}(\mathbb R) \otimes{\cal F}_{t}$-measurable
way, which shows that the functional under discussion,
$(g, \omega) \mapsto\alpha_{t}(g)(\omega)$, is also ${\cal
B}(\mathbb R)
\otimes{\cal F}_{t}$-measurable, being the limit
of the sequence $((g, \omega) \mapsto a^{g,\omega}_{n})_{n=0}^{\infty
}$ of measurable
functions. 

\begin{lemma}
Suppose that for $t \in\{0, \ldots, T-1\}$,
%
\begin{eqnarray}
&& c \dvtx \mathbb R\times\Omega
 \to \mathbb R, \qquad{\cal
B}(\mathbb R)
\otimes{\cal F}\mbox{-measurable such that} \\
&&\quad
a \mapsto c_{t}(a)(\omega) \qquad\mbox{is non-decreasing, continuous,}\label{Cfunct0} \\
 &&\alpha_{t+1} \dvtx \mathbb R\times\Omega
 \to[0, \pi],\qquad
{\cal
B}(\mathbb R) \otimes{\cal F}\mbox{-measurable such that} \label
{alphafunct-10} \\
&&\quad
g \mapsto\alpha_{t+1}(g)(\omega)\qquad
\mbox{is non-decreasing for all }\omega\in\Omega. \label
{alphafunct10}
\end{eqnarray}
Given a regular version $Q_{t}$ of the ${\cal F}_{t}$-conditioned
distribution $Q$, assume that the random variable $\varepsilon_{t+1}$ satisfies
%
\begin{eqnarray}
&&\mathbb R\to[0,\pi],\qquad  x \mapsto\int_{\Omega}\alpha
_{t+1}\bigl(x+\varepsilon_{t+1}(\omega')\bigr)(\omega'){\mathrm Q}_{t}(\mathrm{d}\omega
')(\omega)\nonumber\\[-8pt]\\[-8pt]
&&\phantom{\mathbb R\to[0,\pi],\qquad\;}\quad\mbox{is continuous for each }\omega\in\Omega.\nonumber
\label{measuremu0}
\end{eqnarray}
%
The following assertions then hold:
%
%
\begin{enumerate}[(ii)]
\item[(i)]
there exists a unique
${\cal B}(\mathbb R) \otimes{\cal F}_{t}$-measurable $[0, \pi
]$-valued $\alpha
_{t}$ satisfying
%
\begin{eqnarray} \label{condexp0}
&&\alpha_{t}(g)(\omega)= \mathbb E^{ {\mathrm Q}}_{t}\bigl(\alpha_{t+1}\bigl(g -
c_{t}(\alpha
_{t}(g))+\varepsilon_{t+1}\bigr)\bigr)(\omega)\nonumber\\[-9pt]\\[-9pt]
&&\quad\mbox{for all }g \in\mathbb R,
\mbox{ for almost all }\omega\in\Omega;\nonumber
\end{eqnarray}
\item[(ii)] the mapping $\mathbb R\to[0, \pi]$, $g \mapsto\alpha
_{t}(g)(\omega)$ is non-decreasing and continuous
for all $\omega\in\Omega$.
\end{enumerate}
\label{secondlemma}
\end{lemma}
\begin{pf}
(i) As in the proof of the Lemma \ref{firstlemma}, we obtain the
unique root $\alpha_{t}(g)( \omega)$ of the function
\[
f^{g, \omega}(a):=a
-\int_{\Omega}\alpha_{t+1}\bigl(g-c_{t}(a)(\omega)+\varepsilon
_{t+1}(\omega')\bigr)(\omega'){\mathrm Q}_{t}(\mathrm{d}\omega')(\omega),\qquad  a
\in[0, \pi].
\]
By the
bisection method,
\[
a^{g, \omega}_{n}=\tfrac{1}{2}(\overline{a}^{g, \omega}_{n} +
\underline{a}^{g, \omega}_{n}),\qquad
\begin{array}{lll} \overline{a}^{g, \omega}_{n+1}=a^{g \omega}_{n},
\qquad&\underline{a}^{g, \omega}_{n+1}:=\underline{a}^{g, \omega}_{n}, \qquad&
\mbox{if }f^{g, \omega}(a^{g, \omega}_{n}) \ge0,\vspace*{1pt}\vspace*{2pt}\cr
\overline{a}^{g}_{n+1}=\overline{a}^{g, \omega}_{n}, &\underline
{a}^{g, \omega}_{n+1}:=a^{g}_{n}, & \mbox{if }f^{g, \omega}(a^{g,
\omega}_{n}) < 0,
\end{array}
\]
started at $\underline{a}^{g, \omega}_{0}:=0$, $\overline{a}^{g,
\omega}_{0}:=\pi$.
Since
\[
(g, \omega, a) \mapsto f^{g, \omega}(a)\qquad
\mbox{is }{\cal B}(\mathbb R) \otimes{\cal F}_{t} \otimes{\cal
B}([0, \pi
])\mbox{-measurable},
\]
each bisection point $(g, \omega) \mapsto a^{g, \omega}_{n}$ is
${\cal B}(\mathbb R) \otimes{\cal F}_{t}$-measurable, which shows that
for $n \to\infty$,
the pointwise limit $(g, \omega) \mapsto\alpha_{t}(g, \omega)$ of
the bisection sequence is also ${\cal B}(\mathbb R) \otimes{\cal
F}_{t}$-measurable.
By construction, the equality
\[
\alpha_{t}(g)(\omega)=\int_{\Omega}\alpha_{t+1}\bigl(g - c_{t}(\alpha
_{t}(g)(\omega))(\omega)+\varepsilon_{t+1}(\omega')\bigr)(\omega'){\mathrm
Q}_{t}(\mathrm{d}\omega')(\omega)
\]
holds for all $g \in\mathbb R$ and $\omega\in\Omega$, whose right-hand
side is nothing but
the right-hand side of (\ref{condexp0}) for each $g \in\mathbb R$.

(ii) The proof is obtained from (iii) of the previous lemma by replacing
$\alpha_{1}(\cdot)$, $\alpha_{0}(\cdot)$, $c(\cdot)$ and
$Q(\mathrm{d}\omega')$ by $\alpha_{t+1}(\cdot)(\omega)$,$\alpha_{t}(\cdot
)(\omega)$,
$c_{t}(\cdot)(\omega)$
and $Q_{t}(\mathrm{d}\omega')(\omega)$, respectively, with appropriate
notational adaptations according to the conditioning on $\omega$.
\end{pf}

Finally, we address a solution to (\ref{ReducedProblem}) in the last
point of the following proposition.

\begin{proposition} Consider ${\cal E}_{T}$
under the model assumption (\ref{nopointmassnew}) and the cumulative
abatement volume functions
from (\ref{cfunt})
under (\ref{asscosts}) and (\ref{volumefunction}).
\begin{enumerate}[(ii)]
\item[(i)] Given measure ${\mathrm Q} \sim P$, there exist functionals
%
\begin{equation}
\alpha_{t} \dvtx \mathbb R\times\Omega\to[0, \pi],\qquad {\cal
B}(\mathbb R
) \otimes{\cal F}_{t}\mbox{-measurable for }t=0, \ldots T, \label{measurability}
\end{equation}
which fulfill, for all $g \in\mathbb R$,
%
\begin{eqnarray}
\alpha_{T}(g)&=&\pi1_{[0, \infty[}(g), \label{startrecursion} \\
\alpha_{t}(g)&=& \mathbb E^{ {\mathrm Q}}_{t}\bigl(\alpha_{t+1}\bigl(g -
c_{t}(\alpha
_{t}(g))+\varepsilon_{t+1}\bigr)\bigr), \qquad t=0, \ldots, T-1.
\label{assertion}
\end{eqnarray}
\item[(ii)] There exists a ${\mathrm Q}$-martingale $(A^{*}_{t})_{t=0}^{T}$
which satisfies
%
\begin{equation}\label{terminalvalue}
A^{*}_{T}=\pi1_{\{{\cal E}_{T}-\sum_{t=0}^{T-1}c_{t}(A^{*}_{t}) \ge
0\}}.
\end{equation}
\end{enumerate}
\end{proposition}

\begin{pf}
(i) In this proof, we repeatedly make use of Lemma \ref{secondlemma}.
Let us start
with $t=T-1$ and verify that the assumptions of this
lemma are satisfied. Due to continuity (\ref{nondecrcont}) of the
abatement function, we have
(\ref{Cfunct}).
The properties
(\ref{alphafunct-10}) and (\ref{alphafunct10}) hold for $t=T-1$, by
definition (\ref{startrecursion}).
To show (\ref{measuremu0}), we utilize the specific form of $\alpha_{T}$:
%
\begin{equation} \label{assignment}
x \mapsto\int_{\Omega}\alpha_{T}\bigl(x+\varepsilon_{T}(\omega
')\bigr)(\omega'){\mathrm Q}_{T-1}(\mathrm{d}\omega')(\omega)={\mathrm
Q}_{T-1}(x+\varepsilon_{T}\ge0)(\omega).
\end{equation}
Note that, due to (\ref{nopointmassnew}), there are almost surely no
point masses in the distribution of
\[
\varepsilon_{T}=\sum_{i \in I}\gamma^{i}-\mathbb{E}^{Q}_{T-1}\biggl(\sum_{i \in
I}\gamma^{i}\biggr)
\]
conditioned on ${\cal F}_{T-1}$ (with respect $Q$, since $Q \sim P$).
That is, (\ref{assignment}) is continuous for each
$\omega\in\Omega$, as required in (\ref{measuremu0}). Hence, (i)
of Lemma \ref{secondlemma} yields the functional $\alpha_{T-1}$
satisfying (\ref{assertion}) (with $t=T-1$), as required.
To proceed by induction, we emphasize that (ii)
of Lemma \ref{secondlemma} ensures that $g \mapsto\alpha
_{T-1}(g)(\omega)$ is non-decreasing and continuous for
all $\omega\in\Omega$. That is, for the next step, $t=T-2$, the
assumption (\ref{alphafunct10})
on $\alpha_{T-1}$ is automatically satisfied. Moreover, (\ref
{measuremu0}) now follows, due to the continuity of $g \mapsto\alpha
_{T-1}(g)(\omega)$,
from the pointwise convergence
\[
\lim_{n \to\infty} \alpha_{T-1}\bigl(x_{n}+\varepsilon_{T-1}(\omega
')\bigr)(\omega')=\alpha_{T-1}\bigl(x+\varepsilon_{T-1}(\omega')\bigr)(\omega')
\qquad\mbox{for all }\omega' \in\Omega,
\]
dominated by $\pi$, which holds
for each $(x_{n})_{n \in\mathbb N} \subset\mathbb R$ with $\lim_{n
\to\infty
}x_{n}=x$.
That is, all assumptions of Lemma \ref{secondlemma} are also fulfilled
for $t=T-2$. Proceeding recursively
for $t=T-2, \ldots, 0$, we obtain $(\alpha_{t})_{t=0}^{T}$ with (\ref
{measurability}), (\ref{startrecursion}) and
(\ref{assertion}).

(ii) As suggested by (\ref{r1})--(\ref{r3}), we define, for all
$\omega\in\Omega$,
\begin{eqnarray*}
A^{*}_{t}(\omega) &:=& \alpha_{t}(G_{t}(\omega))(\omega),\\
G_{t+1}(\omega) &:=& G_{t}(\omega)-c_{t}(A^{*}_{t}(\omega))(\omega
)+\varepsilon_{t+1}(\omega), \\
&&\mbox{started at }G_{0}:={\cal E}_{0}.
\end{eqnarray*}
The process $(A^{*}_{t})_{t=0}^{T}$ generated in this way obeys the
terminal condition (\ref{terminalvalue}), in view of
(\ref{startrecursion}).
To show the ${\mathrm Q}$-martingale
property of $(A^{*}_{t})_{t=0}^{T}$,
we calculate, for $t=0, \ldots, T-1$,
\begin{eqnarray*}
\mathbb{E}^{ {\mathrm Q}}_{t}(A^{*}_{t+1})(\omega)&=& \mathbb E^{
{\mathrm Q}}_{t}(\alpha
_{t+1}(G_{t+1}))(\omega)\\
&=& \mathbb E^{ {\mathrm Q}}_{t}\bigl(\alpha_{t+1}\bigl(G_{t}(\omega
)-c_{t}(A^{*}_{t}(\omega
))(\omega)+\varepsilon_{t+1} \bigr)\bigr)(\omega) \\
&=& \mathbb E^{ {\mathrm Q}}_{t}\bigl(\alpha_{t+1}\bigl(G_{t}(\omega
)-c_{t}(\alpha
_{t}(G_{t}(\omega))(\omega))(\omega)+\varepsilon_{t+1}
\bigr)\bigr)(\omega)\\
&=&\alpha_{t}(G_{t+1}(\omega) )(\omega) =A^{*}_{t}(\omega)
\end{eqnarray*}
for almost
all $\omega\in\Omega$,
where the penultimate equality follows from (\ref{assertion}).
\end{pf}

\section{Applications}
Let us elaborate on the computational feasibility of our reduced-form modeling.
For illustrative
purposes, we focus on the simplest case of martingales with independent
increments and deterministic abatement functions.
We assume that:
%
\begin{eqnarray} \label{independence}
&&\varepsilon_{t+1}\mbox{ and }{\cal F}_{t}\mbox{ are independent under }Q
\mbox{ for all }t=0, \ldots, T-1; \\[-2pt]
&& c_{t} \dvtx[0, \infty[ \to\mathbb R\mbox{ is deterministic and time
constant }(c_{t}=c)_{t=0}^{T-1}.
\end{eqnarray}
Under these assumptions, the randomness enters the allowance price
through the present up-to-day emissions only. More precisely,
(\ref{independence}) ensures that
%
\begin{equation}
\omega\mapsto\alpha_{t}(g)(\omega) =\alpha_{t}(g) \qquad\mbox{is
constant on }\Omega. \label{detprop}
\end{equation}
Let us verify this assertion.
For $t=T$, (\ref{detprop}) holds, by definition (\ref{startrecursion}).
For $t=T-1, \ldots, 1$, we proceed inductively as follows: by
construction, $\alpha_{t}(g)(\omega)$ is the unique solution $a$ to
%
\begin{eqnarray}
a&=&
\int_{\Omega}\alpha_{t+1}\bigl(g-c_{t}(a)+\varepsilon_{t+1}(\omega
')\bigr)(\omega'){\mathrm Q}_{t}(\mathrm{d}\omega')(\omega) \nonumber\\
&=&
\int_{\Omega}\alpha_{t+1}\bigl(g-c_{t}(a)+\varepsilon_{t+1}(\omega
')\bigr){\mathrm Q}(\mathrm{d}\omega'), \label{notdependent}
\end{eqnarray}
where, in the last equality, we have utilized the fact that $Q_{t}=Q$,
due to the independence (\ref{independence}),
and the fact that $\alpha_{t+1}(g)$ does not depend on $\omega$, by the
the induction assumption.
Obviously, the fixed point $\alpha_{t}(g)(\omega):=a$ from (\ref
{notdependent}) also does not depend on $\omega$.

For numerical calculation, we rely on the one-dimensional least-squares
Monte Carlo method,
which is applicable in our case of martingales with independent
increments. Although this setting is relatively
restrictive, it covers a sufficiently rich class of martingales. For
instance, important cases
of information shocks leading to allowance price jumps can be easily
addressed under this approach when
$({\cal E}_{t})_{t=0}^{T}$ is modeled as an
appropriately sampled, centered Poisson process. In this case, fixed
point equations can be
treated analytically. We do not follow this path in favor of numerical methods,
which deserve particular attention due to the complexity of emissions
markets. In particular,
extensions of Monte Carlo methods to the multidimensional setting (see
\cite{Stentoft2004})
seem to be appropriate. A preliminary analysis shows that assuming the
existence of a global Markovian
state process allows independence to be weakened to conditional
independence, which leads to
multidimensional Monte Carlo, in the sense of \cite{Stentoft2004},
since the state process gives additional dimensions.

We now focus on computational aspects.
From (\ref{detprop}), it follows that $\alpha_{t}(G_{t})$ is a
$\sigma(G_{t})$-measurable random variable.
Thus, in the equality (\ref{assertion}), the
condition ${\cal F}_{t}$ can be replaced by the condition $\sigma(G_{t})$:
%
\begin{equation} \label{fixpoint}
\alpha_{t}(G_{t})= \mathbb E^{Q}\bigl(\alpha_{t+1}\bigl(G_{t}-c_{t}(\alpha
_{t}(G_{t}))+\varepsilon_{t+1}\bigr)  \vert \sigma(G_{t})\bigr).
\end{equation}
We shall treat this relation as a fixed point equation for the
Borel-measurable function $\alpha_{t}$ and attempt to obtain a solution
in the limit $\alpha_{t}=\lim_{n \to\infty}\alpha^{n}_{t}$ of iterations
%
\begin{equation}
\alpha^{n+1}_{t}(G_{t})=\mathbb{E}^{Q}\bigl(\alpha_{t+1}\bigl(G_{t}-c_{t}(\alpha
^{n}_{t}(G_{t}))+\varepsilon_{t+1}\bigr)  |  \sigma(G_{t})\bigr),\qquad
 n \in\mathbb N,
\label{indeeddefine}
\end{equation}
started at $\alpha^{0}_{t}=\alpha_{t+1}$. (Note that, given $\alpha
_{t+1}$ and $\alpha^{n}_{t}$, the equation (\ref{indeeddefine})
indeed defines
a Borel function $\alpha^{n+1}_{t}$ by the factorization of the
$\sigma(G_{t})$-measurable random variable on the right-hand side of
(\ref{indeeddefine}).)
For numerical calculation of conditional expectations, we suggest using
the least-squares Monte Carlo method.

To explain the principle of the least-squares
Monte Carlo approach (see \cite{Longstaff} and \cite{Stentoft2004})
in more detail, we abstract from the concrete situation (\ref
{indeeddefine}) and consider
\[
\varphi(G)= \mathbb E^{Q}(\phi(G, \varepsilon)  |  \sigma(G) ),
\]
where $\varepsilon, G$ are $\mathbb R$-valued and independent with
respect to $Q$ and $\phi$ is a bounded Borel function on $\mathbb R^{2}$.
Under these assumptions, the function
$\varphi$ is obtained as
$
\varphi(g)=\int_{\mathbb R}\phi(g, e)Q^{\varepsilon}(\mathrm{d}e)
$ for $Q^{G}$-almost all $g \in\mathbb R$,
where $Q^{\varepsilon}, Q^{G}$ are image measures of $Q$ under
$\varepsilon$ and $G$, respectively.
An equivalent condition defining $\varphi$ is the orthogonality
%
\begin{eqnarray}\label{strong}
&&\mbox{determine }\varphi\in L^{2}(\mathbb R, \mu)\mbox{ such that for all }
\psi\in\Psi,\nonumber\\[-11pt]\\[-11pt]
&&\int_{\mathbb R^{2}}\bigl(\varphi(g)-\phi(g, e)\bigr) \psi(g) (Q^{\varepsilon
}\otimes\mu)(\mathrm{d}e,\mathrm{d}g)=0,\nonumber
\end{eqnarray}
where $\mu$ is a measure which is equivalent to $ Q^{G}$ and $\Psi$
stands for a set of functions
which are square-integrable with respect to $\mu$, whose
linear space is dense in $L^{2}(\mathbb R, \mu)$. The idea of the
least-squares Monte Carlo method is to relax,
for computational tractability, the principle (\ref{strong})~to
%
\begin{eqnarray}\label{weak}
&&\mbox{determine }\varphi\in\operatorname{lin}\Psi\mbox{ such\ that\ for\ all }\psi
\in\Psi,\nonumber\\[-11pt]\\[-11pt]
&&\sum_{k=1}^{K}\bigl(\varphi(g_{k})-\phi(g_{k}, e_{k})\bigr) \psi(g_{k})=0,\nonumber
\end{eqnarray}
with a finite set of basis functions
\[
\Psi=\{\psi_{j}  \dvtx j=1, \ldots, J\}
\]
and an appropriate sample
\[
S:= ( e_{k},g_{k})_{k=1}^{K} \subset\mathbb R^{2},
\]
chosen such that
the combination
$
\frac{1}{K}\sum_{k=1}^{K}\delta_{(e_{k}, g_{k})}
$
of the Dirac measures approximates the distribution $Q^{\varepsilon
}\otimes\mu$
(for instance, $S$ being realizations of $K \in\mathbb N$
independent $Q^{\varepsilon}\otimes\mu$-distributed random
variables). The solution
to the weakened problem (\ref{weak}) is given in terms of
%
\begin{eqnarray}
\label{realizSample}
&&\mbox{realizations }\phi(S)=(\phi(e_{k},g_{k}))_{k=1}^{K}\mbox{  of  }\phi
\mbox{ on the sample $S$}, \nonumber\\[-10pt]\\[-10pt]
&&\mbox{realizations }M= ( \psi_{j}(g_{k})
)_{k=1,j=1}^{K,J}\mbox{ of basis functions on $S$},\nonumber
\end{eqnarray}
as follows:
\begin{eqnarray*}
&&\mbox{if }q=(q_{j})_{j=1}^{J} \mbox{ fulfills  }M^{\top}Mq=M^{\top}\phi
(S),\\[-4pt]
&&\mbox{then (\protect\ref{weak}) is solved by }
\varphi= \sum_{j=1}^{J}q_{j}\psi_{j}.
\end{eqnarray*}
%
\begin{figure}[b]

\includegraphics{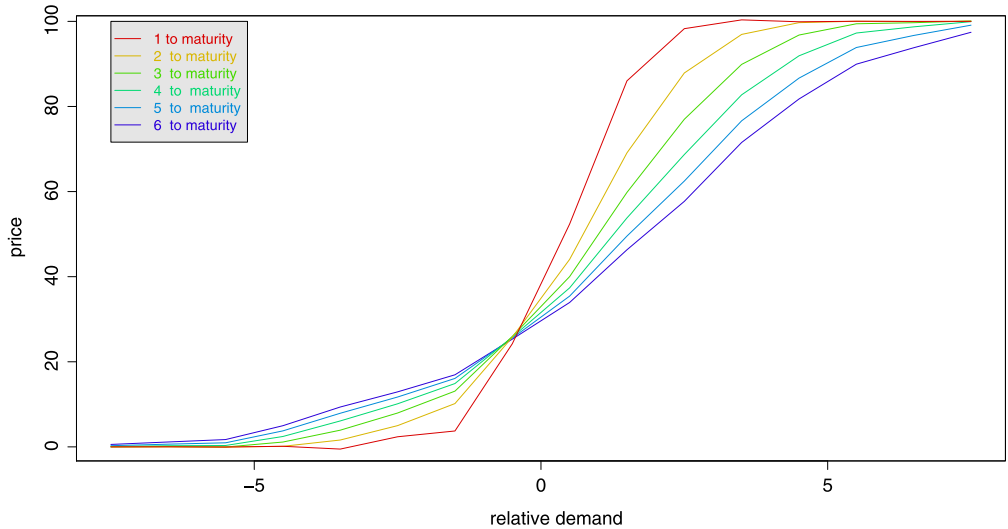}

\caption{The functions $\alpha_{t}$ for $t=T-1, \ldots, T-6$, from the
least-squares Monte Carlo method.}\label{alphafigure}\vspace*{-5pt}
\end{figure}
We now formulate an algorithm for the approximate calculation of (\ref
{fixpoint}) in which the conditional
expectation is replaced by the least-squares Monte Carlo projection. To
ease notation, let us suppose
that $(\varepsilon_{t})_{t=1}^{T}$ are identically distributed (in
addition to their independence (\ref{independence})).

\begin{APMCM*}
\begin{enumerate}
\item\textit{Initialization.}
Given sample $S=(e_{k}, g_{k})_{k=1}^{K} \subset\mathbb R^{2}$
describing the distribution of $Q^{\varepsilon_{1}} \otimes\mu$
and a set of basis functions $\Psi=(\psi_{i})_{j=1}^{J}$ on $\mathbb R$,
define $M$ as in (\ref{realizSample}).
Set $\alpha_{T}(g)=1_{[0, \infty[}(g)$ for all $g \in\mathbb R$ and
proceed in the next step with $t:=T-1$.
\item\textit{Iteration.} Define $\alpha^{0}_{t}=\alpha_{t}$ and proceed
in the next step with $n:=0$.
\begin{enumerate}[(2b)]
\item[(2a)] Calculate $\phi^{n+1}(S):=(\alpha
_{t+1}(g_{k}-c_{t}(\alpha^{n}_{t}(g_{k}))+e_{k}))_{k=1}^{K}$.
\item[(2b)] Determine a solution $q^{n+1} \in\mathbb R^{J}$ to
$
M^{\top}Mq^{n+1}=M^{\top}\phi^{n+1}(S)
$.
\item[(2c)] Define $\alpha^{n+1}_{t}:=\sum_{j=1}^{J}q^{n+1}_{j}\psi
_{j}$.
\item[(2d)] If $\max_{k=1}^{K}|\alpha^{n+1}_{t}(g_{k})-\alpha
^{n}_{t}(g_{k})|\ge\varepsilon$, then
put $n:=n+1$ and continue with step
2a).

\noindent
If $\max_{k=1}^{K}|\alpha^{n+1}_{t}(g_{k})-\alpha
^{n}_{t}(g_{k})|<\varepsilon$, then set $t:=t-1$.
If $t>1$, go to step 2, otherwise finish.
\end{enumerate}
\end{enumerate}
\end{APMCM*}

\begin{example*}
To illustrate allowance price calculation via the Monte Carlo method,
we consider the following
numerical example. Suppose that the penalty
 is set at
$\pi=100$
and that the martingale increments
$(\varepsilon_{t})_{t=1}^{T}$ are independent, identically Normally
distributed. Note that by an appropriate choice of the
emission measurement scale, the standard deviation can always be
normalized, thus we have assumed that
each $\varepsilon_{t}$ is ${\cal N}(0.5, 1)$-distributed. Further,
consider the
basis consisting of piecewise linear hut functions
\[
\psi_{j} \dvtx \mathbb R\to[0, 1],\qquad  x \mapsto(1-|z_{j}-x|/h)^{+}\qquad
\mbox{for }x \in\mathbb R,\ j=1, \ldots, J,
\]
where the peaks $z_{1}=-(J-1)*h/2, \ldots, z_{J}=(J-1)*h/2$ are chosen
to be equidistant with the distance
$h>0$. For numerical illustration,
we set $J=16$ and $h:=1$. Further, the sample $S=(e_{k},
g_{k})_{k=1}^{K}$ for the Monte Carlo method is generated
with $K=1000$ outcomes. For $(e_{k})_{k=1}^{K}$, we followed a natural
choice, taking realizations of
$K$ independent ${\cal N}(0.5, 1)$-distributed random variables.
However, since
the distribution of $G_{t}$ is not known in advance, an appropriate
candidate for
$\mu$ seems to be the uniform distribution concentrated on the
interval which is relevant for calculations.
That is, the outcomes $(g_{k})_{k=1}^{K}$ are constructed by
equidistant sampling of $[z_{1}, z_{J}]$,
ranging from $g_{1}=z_{1}=-7.5$ to $g_{K}=z_{J}=7.5$. For
the cumulative volume function $c \dvtx \mathbb R\to\mathbb R, a \mapsto
0.1\sqrt
{(a)^{+}}$, we observed a fast and stable convergence
which gave a reasonable outcome within a few iterations. The resulting
functions $(\alpha_{t})_{t=T-1}^{T-6}$ are depicted in Figure~\ref
{alphafigure}.
\end{example*}

Let us outline a valuation procedure for a European call on emission
allowance price.

\begin{VECMCM*}
\begin{enumerate}
\item Given basis functions $\Psi=(\psi_{j})_{j=1}^{J}$ and a sample
$S=(e_{k}, g_{k})_{k=1}^{K} \subset\mathbb R^{2}$
which approximates $Q^{\varepsilon_{1}} \otimes\mu$,
determine $(\alpha_{t})_{t=T}^{0}$ in terms of basis coefficients
using the above
 least-squares Monte Carlo
algorithm.
\item Given maturity time $\tau\in\{1, \ldots, T\}$ of the European
call, determine its pay-off
$f^{\tau}_{\tau}:=(\alpha_{\tau}-K)^{+}$.
Calculate least-squares projections, recursively processing for $u=\tau
, \ldots, t $ as follows:
\begin{enumerate}[(b)]
\item[(a)] put $\phi(S)=(f^{\tau}_{u}(g_{k}-c_{u}(\alpha
_{u}(g_{k}))+e_{k}))_{k=1}^{K}$;
\item[(b)] obtain $q$ as solution to $M^{\top}Mq=M\phi$;
\item[(c)] set $f^{\tau}_{u-1}=\sum_{j=1}^{J}q_{j}\psi_{j}$;
\item[(d)] if $u-1=t$, then finish, otherwise set $u:=u-1$ and return
to (a).
\end{enumerate}
\item Given recent allowance price $a$, calculate the state variable
$g$ as solution to
$a=\alpha_{t}(g)$.
\item Plug in the state variable $g$ and into function $f^{\tau}(t,
\cdot)$ to obtain the price of the European call as
as $f^{\tau}(t, g)$.
\end{enumerate}
\end{VECMCM*}

Let us conclude this section by sketching core ideas on continuous-time
modeling.
Our analysis shows that the risk-neutral allowance price evolution
$(A_{t})_{t=0}^{T}$ must be described by
a martingale whose
terminal value is digital and depends on the intermediate values (see
(\ref{ReducedProblemCont})).
Suppose that the compliance period
is given by an interval $[0, T]$, such that all relevant
random evolutions are described by adapted stochastic processes on
%
\begin{eqnarray} \label{equipped}
&&\mbox{filtered
probability space }(\Omega, {\cal F}, P, ({\cal F}_{t})_{t \in[0,
T]})\nonumber\\[-11pt]\\[-11pt]
&&\mbox{equipped with probability measure }Q \sim P,\nonumber
\end{eqnarray}
where $Q$ represents the spot martingale measure.
Given a random variable ${\cal E}_{T}$ and appropriate non-decreasing and
continuous abatement functions $c_{t}: \mathbb R_{+} \times\Omega\to
\mathbb R_{+}$
indexed by $t \in[0, T]$, we follow
an analogy to discrete time and consider
solutions $(A_{t})_{t \in[0, T]}$ to
%
\begin{equation} \label{charact1}
A_{t}=\pi E^{Q}\bigl(1_{\{{\cal E}_{T}-\int_{0}^{T}c_{t}(A_{s})\,\mathrm{d}s
\ge0\}}
 |
 {\cal F}_{t}\bigr),\qquad  t \in[0, T].
\end{equation}
%
Our results from the discrete-time setting suggest that if
%
\begin{eqnarray}
\cases{\mbox{the increments of the martingale }
({\cal E}_{t}=E^{Q}({\cal E}_{T} |  {\cal F}_{t}))_{t \in[0, T]}\mbox{ are}
\cr
\mbox{independent and the abatement functions }c_{t} \dvtx \mathbb R_{+}
\times
\to\mathbb R_{+} \cr
\mbox{are deterministic and time constant }(c_{t}=c)_{t \in[0, T]},
}\label{emart}
\end{eqnarray}
then a solution to (\ref{charact1}) should be expected in
the functional form
\[
A_{t}=\alpha(t, G_{t}),\qquad  t \in[0, T],
\]
with an appropriate deterministic function
%
\begin{equation}
\alpha \dvtx  [0, T]\times\mathbb R\mapsto\mathbb R,\qquad (t, g)
\mapsto
\alpha(t,g),
\end{equation}
and a state process $(G_{t})_{t \in[0, T]}$ given by
\[
G_{t}={\cal E}_{t}-\int_{0}^{t}c_{s}(A_{s})\,\mathrm{d}s,\qquad  t \in[0, T].
\]
To illustrate how such an approach allows one to guess a solution,
assume that
%
\begin{eqnarray} \label{supports}
&&(\Omega, {\cal F}, P, ({\cal F}_{t})_{t \in[0, T]})
\mbox{ supports the process }(W_{t}, {\cal F}_{t})_{t \in[0, T]} \nonumber\\[-11pt]\\[-11pt]
&&\mbox{of Brownian motion with respect to } Q\sim P. \nonumber
\end{eqnarray}
Furthermore, we respond to (\ref{emart}), supposing that
\begin{eqnarray*}
&&\mathrm{d}{\cal E}_{t}=\sigma_{t}\,\mathrm{d}W_{t} \mbox{with pre-specified
deterministic }(\sigma_{t})_{t \in[0, T]}\mbox{ and} \\[-4pt]
&&\mbox{continuous and non-decreasing abatement functions }
(c_{t}=c)_{t \in[0, T]}.
\end{eqnarray*}
To ensure the martingale property of
$
(A_{t}=\alpha(t, G_{t}))_{t \in[0, T]}
$,
apply the It\^o formula
\begin{eqnarray*}
\mathrm{d}A_{t}
&=&\mathrm{d}\alpha(t, G_{t})
=\partial_{(1, 0)}\alpha(t, G_{t})\,\mathrm{d}t+ \partial_{(0, 1)}\alpha(t,
G_{t})\,\mathrm{d} G_{t}
+ \tfrac{1}{2}\partial_{(0, 2)}\alpha(t, G_{t})\,\mathrm{d}[G]_{t}\\
&=& \partial_{(1, 0)}\alpha(t, G_{t})\,\mathrm{d}t-
\partial_{(0, 1)}\alpha(t, G_{t})c(\alpha(t,G_{t}))\,\mathrm{d}t+\tfrac
{1}{2}\partial_{(0, 2)}\alpha(t, G_{t})\sigma^{2}_{t}\,\mathrm{d}t \\
&&{}
+ \partial_{(0, 1)}\alpha(t, G_{t}) \sigma_{t}\,\mathrm{d}W_{t}
\end{eqnarray*}
and claim the function $\alpha$ as a solution on $[0, T[ \times
\mathbb R$ to
%
\begin{equation} \label{pde}
\partial_{(1, 0)}\alpha(t, g)-
\partial_{(0, 1)}\alpha(t, g)c(\alpha(t,g))+\tfrac{1}{2}\partial
_{(0, 2)}\alpha(t, g)\sigma^{2}_{t}=0
\end{equation}
with boundary condition
%
\begin{equation} \label{bc}
\alpha(T, g)=\pi1_{[0, \infty[}(g)\qquad\mbox{for all }g \in\mathbb
R,
\end{equation}
justified by the
digital terminal allowance price.
Having obtained $\alpha$ in this way, we construct the state process
as the solution to the stochastic differential equation
%
\begin{equation}
\mathrm{d}G_{t}=\mathrm{d}{\cal E}_{t}-c(\alpha(t, G_{t}))\,\mathrm{d}t, \qquad G_{0}={\cal
E}_{0},
\label{sde}
\end{equation}
and then determine
%
\begin{equation}
A_{t}:=\alpha(t, G_{t}),\qquad  t \in[0, T]. \label{introd}
\end{equation}
Finally, this process must be verified in order to solve (\ref{charact1}).

\section{Conclusion}
This article explains the logical principles underlying risk-neutral
modeling of emission certificate
price evolution.
We show that within a
realistic situation of risk-averse market players,
there is no connection between social optimality and market
equilibrium, but there is a
useful feedback relation characterizing
risk-neutral allowance price dynamics.
Expressing this result in terms of fixed point equations on the level
of martingales,
we address the existence of its solution and elaborate on its
algorithmic tractability.
Furthermore, we suggest an extension of these concepts to continuous time
and show that promising results can be obtained using
diffusion processes. Here, emission allowances and their options
can be described in terms of standard partial differential equations.
Although option pricing in this framework seems to be appealing,
we believe that it is not
superior to our Monte Carlo method since the latter can be used in high
dimensions and, more importantly,
in the presence of jumps
in the martingale $({\cal E}_{t})_{t=0}^{T}$. This is particularly
important to describe price shocks,
which may result from possible discontinuities in the information flow.

\section*{Acknowledgements}
The authors would like to thank the referees
for insightful remarks
and comments which helped us to improve this work.

\printhistory


\begin{thebibliography}{99}

\bibitem{BenzTrueck}
Benz, E. and Trueck, S. (2008).
Modeling the price dynamics of co2 emission allowances.
\textit{Energy Economics} \textbf{31} 4--15.

\bibitem{CarmonaFehrHinz}
Carmona, R., Fehr, F. and Hinz, J. (2009).
Optimal stochastic control and carbon price formation.
\textit{SIAM J.~Control Optim.}
\textbf{48} 2168--2190.
\MR{2520324}

\bibitem{CarmonaFehrHinzPorchet}
Carmona, R., Fehr, F., Hinz, J. and Porchet, A. (2010).
Market designs for emissions trading schemes.
  \textit{SIAM Review}. To appear.

\bibitem{CarmonaHinz}
Carmona, R. and Hinz, J. (2009).
Risk neutral modeling of emission allowance prices and option
valuation.
Technical report, Princeton Univ.

\bibitem{CetinVerschuere}
Cetin, U. and Verschuere, M. (2009).
Pricing and hedging in carbon emissions markets.
\textit{Int. J. Theor. Appl. Finance}
\textbf{12}  949--967.

\bibitem{ChesneyTaschini}
Chesney, M. and Taschini, L. (2008).
The endogenous price dynamics of the emission allowances: An
application to co2 option pricing.
Technical report.

\bibitem{Cronshaw}
Cronshaw, M. and Kruse, J.B. (1996).
Regulated firms in pollution permit markets with banking.
\textit{Journal of Regulatory Economics} \textbf{9} 179--189.

\bibitem{Dales}
Dales, J.H. (1968).
\textit{Pollution, Property and Prices}. Toronto:
Univ. Toronto Press.

\bibitem{DaskalakisPsychoyiosMarkellos}
Daskalakis, G., Psychoyios, D. and Markellos, R.N. (2009).
Modeling co2 emission allowance prices and derivatives: Evidence from
the European Trading Scheme. \textit{Journal of Banking and Finance} \textbf{33}
1230--1241.

\bibitem{Kabanov}
Kabanov, Y.M. and Stricker, C. (2001).
A teachers' note on no-arbitrage criteria.
\textit{Lecture Notes Math.} \textbf{1775} 149--152.
\MR{1837282}

\bibitem{Leiby}
Leiby, P. and Rubin, J. (2001).
Intertemporal permit trading for the control of greenhouse gas
emissions.
\textit{Environmental and Resource Economics} \textbf{19} 229--256.

\bibitem{Longstaff}
Longstaff, F. and Schwartz, E. (2001).
Valuing american options by simulation: A simple least-squares
approach.
\textit{Review of Financial Studies} \textbf{14} 113--147.

\bibitem{Maeda}
Maeda, A. (2004).
Impact of banking and forward contracts on tradable permit markets.
\textit{Environmental Economics and Policy Studies} \textbf{6} 81--102.

\bibitem{Montgomery}
Montgomer, W.D. (1972).
Markets in licenses and efficient pollution control programs.
\textit{Journal of Econom. Theory} \textbf{5} 395--418.
\MR{0443849}

\bibitem{PaolellaTaschini}
Paolella, M.S. and Taschini, L. (2008).
An econometric analysis of emissions trading allowances.
\textit{Journal of Banking and Finance}
\textbf{32} 2022--2032.

\bibitem{Rubin}
Rubin, J. (1996).
A model of intertemporal emission trading, banking and borrowing.
\textit{Journal of Environmental Economics and Management}
\textbf{31} 269--286.

\bibitem{Schennach}
Schennach, S.M. (2000).
The economics of pollution permit banking in the context of title iv
of the 1990 clean air act amendments.
\textit{Journal of Environmental Economics and Management}
\textbf{40} 189--210.

\bibitem{Seifert}
Seifert, J. Uhrig-Homburg, M. and Wagner, M. (2008).
Dynamic behavior of co2 spot prices.
\textit{Journal of Environmental Economics and Management} \textbf{56}
180--194.

\bibitem{sijm}
Sijm, J., Neuhoff, K. and Chen, Y. (2006).
Co2 cost pass-through and windfall profits in the power.
\textit{Climate Policy} \textbf{6} 49--72.

\bibitem{Stentoft2004}
Stentoft, L. (2004).
Convergence of the least squares monte carlo approach to american
option valuation.
\textit{Management Science} \textbf{50} 576--611.

\bibitem{StevensRose}
Stevens, B. and Rose, A. (2002).
A dynamic analysis of the marketable permits approach to global
warming policy: A comparison of spatial and temporal flexibility.
\textit{Journal of Environmental Economics and Management} \textbf{44}
45--69.

\bibitem{Tietenberg2}
Tietenberg, T. (1985).
\textit{Emissions Trading: An Exercise in Reforming Pollution Policy}.
Boston: Resources for the Future.

\bibitem{UhrigHomburgWagner1}
Uhrig-Homburg, M. and Wagner, M. (2008).
Derivatives instruments in the {E}{U} emissions trading scheme, and
early market perspective.
\textit{Energy and Environment} \textbf{19} 635--655.

\bibitem{UhrigHomburgWagner2}
Uhrig-Homburg, M. and Wagner, M. (2009).
Futures price dynamics of c{O}2 emissions certificates: An empirical
analysis.
\textit{Journal of Derivatives}  \textbf{17}  73--88.

\bibitem{Wagner}
Wagner, M. (2006).
co2-Emissionszertifikate, Preismodellierung und
Derivatebewertung.
Ph.D. thesis, Universit\"at Karlsruhe.

\end{thebibliography}
\end{document}